\documentclass{amsart}
\usepackage{amssymb}
\usepackage{color} 

\newtheorem{theorem}{Theorem}
\newtheorem{lemma}{Lemma}

\begin{document}
\title{Graded codimensions of Lie superalgebra $b(2)$}

\author[D. Repov\v s, and M. Zaicev]
{Du\v san Repov\v s and Mikhail Zaicev}

\address{Du\v san Repov\v s \\Faculty of Education, and
Faculty of  Mathematics and Physics, University of Ljubljana,
P.~O.~B. 2964, Ljubljana, 1001, Slovenia}
\email{dusan.repovs@guest.arnes.si}

\address{Mikhail Zaicev \\Department of Algebra\\ Faculty of Mathematics and
Mechanics\\  Moscow State University \\ Moscow,119992, Russia}
\email{zaicevmv@mail.ru}

\thanks{The first author was supported by the Slovenian Research Agency
grants P1-0292-0101, J1-6721-0101 and J1-5435-0101.
The second author was partially supported by RFBR grant No 13-01-00234a. We thank the referee for comments}



\begin{abstract}
We study asymptotic behaviour of graded codimensions of Lie superalgebra $b(2)$.
We prove that graded PI-exponent exists and is equal to $3+2\sqrt 3$.
\end{abstract}


\maketitle

\section{Introduction}

We consider finite dimensional Lie superalgebras over a field of characteristic zero
and study their $\mathbb{Z}_2$-graded identities. We pay main attention  to
numerical invariants of identities, in particular, to graded codimensions and their 
asymptotic behaviour.

It is well-known that in case $\dim L<\infty$ both graded and ordinary codimensions
are exponentially bounded (\cite{BD-LAA}). One of the more important questions of the
theory of numerical invariants of polynomial identities is: does a (graded)
PI-exponent exist?

There are many papers where the existence of PI-exponent is proved for different 
classes of algebras. For example, if $A$ is an associative PI-algebra or a finite 
dimensional Lie, Jordan or alternative algebra then its PI-exponent exists and is
a non-negative integer (see \cite{GShZ}, \cite{GZ1}, \cite{GZ2}, \cite{Z0}). 
The existence of
PI-exponent for any finite dimensional simple algebra was proved in \cite{GZ-LMS}.
It is not difficult to show that if PI-exponent of $A$ exists then it is less than 
or equal to $d$
provided that $d=\dim A<\infty$ (see for example \cite{BD-LAA}). In many important
classes of algebras over an algebraically closed field (associative, Lie, Jordan,
alternative) the equality $exp(A)=\dim A$ is equivalent
to simplicity of $A$ (\cite{GShZ}, \cite{GZ1}, \cite{Z0}). Recently \cite{GZ-LMS}
it was shown that $exp(L)<\dim L$ for any finite dimensional simple Lie
superalgebra $L$ of the type $b(t), t\ge 3$, (in the notation of \cite{Sch}). The
existence of PI-exponent and similar inequality $exp(L)<\dim L$ for $b(2)$
 was also proved in \cite{GZ-LMS}  although it is not simple superalgebra.

Graded codimensions of Lie superalgebras were studied much less. In particular, it 
is still unknown if $exp^{gr}(L)$ exists even when $L$ is a finite dimensional simple 
Lie superalgebra. In recent paper \cite{RZ-ART} it was proved that an upper graded
PI-exponent of Lie superalgebra $b(t), t\ge 2$, is less than or equal to 
$t^2-1+t\sqrt{t^2-1}$. In particular, this gives an upper bound for ordinary PI-exponent
of $b(t)$. In the present paper we prove the existence of graded PI-exponent of $b(2)$.
We also prove that $exp^{gr}(b(2))=3+2\sqrt{3}$. 
Note that it was recently announced that an ordinary PI-exponent does not exist in
 general non-associative case (see \cite{Z}).

\section{Main constructions and definitions}

Let $L$ be a Lie superalgebra over a field $F$ of characteristic zero, that is
$L=L_0\oplus L_1$ is a non-associative $\mathbb{Z}_2$-graded algebra  satisfying two 
identical relations
$$
xy+(-1)^{|x||y|}yx=0,
$$

$$
x(yz)=(xy)z+(-1)^{|x||y|}y(xz)=0
$$
where $x,y,z\in L_0\cup L_1$ and $|x|=0$ if $x\in L_0$ while $|x|=1$ if $x\in L_1$.

Elements from
 $L_0\cup L_1$ are called {\it homogeneous} and we say that $x$ is 
 {\it even} if
$x\in L_0$ or $x$ is {\it odd} if $x\in L_1$.

Denote by $\mathcal{L}(X,Y)$ a free Lie superalgebra with infinite sets of even
generators $X$ and odd generators $Y$. A polynomial 
$f=f(x_1,\ldots, x_m,y_1,\ldots,y_n) 
\in \mathcal{L}(X,Y)$ is said to be a {\it graded 
identity} of Lie superalgebra 
$L=L_0\oplus L_1$ if $f(a_1,\ldots, a_m,b_1,\ldots,b_n)$ $=0$ 
whenever $a_1,\ldots, a_m\in 
L_0, b_1,\ldots,b_n\in L_1$.

Given positive integers $0\le k\le n$, denote by $P_{k,n-k}$ the subspace of all 
multilinear polynomials  $f=f(x_1,\ldots, x_k,y_1,\ldots,y_{n-k})\in \mathcal{L}(X,Y)$
of degree $k$ on even variables and of degree $n-k$ on odd variables. Denote by
$Id^{gr}(L)$ the ideal of $\mathcal{L}(X,Y)$ of all graded identities of $L$. Then
$P_{k,n-k}\cap Id^{gr}(L)$ is the subspace of all multilinear graded identities of $L$
of total degree $n$ depending on $k$ even variables and $n-k$ odd variables. Also denote
by $P_{k,n-k}(L)$ the quotient
$$
P_{k,n-k}(L)=\frac{P_{k,n-k}}{P_{k,n-k}\cap Id^{gr}(L)}.
$$
Then the graded $(k,n-k)$-codimension of $L$ is
$$
c_{k,n-k}(L)=\dim P_{k,n-k}(L)
$$
and the total graded codimension of $L$ is
$$
c_n^{gr}(L)=\sum_{k=0}^n{n\choose k} c_{k,n-k}(L).
$$

It is known (see \cite{BD-LAA}) that if $\dim L <\infty$ then the sequence 
$\{c_n^{gr}(L)\}_{n\ge 1}$ is  exponentially bounded and one can 
consider the related 
sequence $\sqrt[n]{c_n^{gr}(L)}$. The latter sequence
 has the lower and upper limits 
$$
\underline{exp}^{gr}(L)=\liminf_{n\to\infty} \sqrt[n]{c_n^{gr}(L)},
\qquad
\overline{exp}^{gr}(L)=\limsup_{n\to\infty} \sqrt[n]{c_n^{gr}(L)}
$$
called the {\it lower}  and {\it upper graded} PI-{\it exponents} of $L$, 
respectively. If an ordinary limit exists, it is called an 
{\it (ordinary) graded} PI-{\it exponent} of $L$,
$$
exp^{gr}(L)=\lim_{n\to\infty} \sqrt[n]{c_n^{gr}(L)}.
$$

Symmetric groups and their representations play an important role in the theory of 
codimensions. One can find all details concerning application of representation theory of
symmetric groups to study of polynomial identities  in \cite{B}, 
\cite{Dren}, \cite{GZbook}. In case of graded identities one can consider 
$(S_k\times S_{n-k})$-action on multilinear graded polynomials. Namely, the subspace
$P_{k,n-k}\subseteq\mathcal{L}(X,Y)$ has a natural structure of $(S_k\times S_{n-k})$-
module where $S_k$ acts on even variables $x_1,\ldots,x_k$ while $S_{n-k}$ acts
on odd variables $y_1,\ldots,y_{n-k}$. Clearly, $P_{k,n-k}\cap Id^{gr}(L)$ is the
submodule under this action and we get an induced $S_k\times S_{n-k}$-action on 
$P_{k,n-k}(L)$. If $G$ is a subgroup of $S_k\times S_{n-k}$ then $G$ also acts 
naturally on $P_{k,n-k}(L)$. In particular,
\begin{equation}\label{eq0}
c_{k,n-k}(L)\ge \dim M
\end{equation}
for any subgroup $G\subseteq S_k\times S_{n-k}$ and for any $G$-submodule $M$ of
 $P_{k,n-k}(L)$. We will use the relation (\ref{eq0}) for getting a lower bound 
 of $c_{k,n-k}(L)$.

\section{Graded PI-exponent of $b(2)$}

Recall the construction of Lie superalgebra $L=L_0\oplus L_1=b(2)$. Even component $L_0$
consists of all $4\times 4$ matrices of the type
$$
L_0 = \left\{ \left(
           \begin{array}{cc}
             A & 0 \\
             0 & -A^t \\
           \end{array}
         \right)
 \mid A\in M_2(F), tr(A)=0  \right\},
$$
where $A$ is $2\times 2$ traceless matrix and $t:A\to A^t$ is the usual transpose involution.

Odd component $L_1$ consists of $4\times 4$ matrices
$$
L_1 = \left\{ \left(
           \begin{array}{cc}
             0 & B \\
             C & 0 \\
           \end{array}
         \right)
 \mid{\rm where}\quad B^t=B, C^t=-C\,\quad B,C\in M_2(F)  \right \}.
$$
Also denote
$$
L^{-}_1 = \left\{ \left(
           \begin{array}{cc}
             0 & 0 \\
             C & 0 \\
           \end{array}
         \right) \mid C^t=-C\in M_2(F)  \right\},
$$
$$
L^{+}_1 = \left\{ \left(
           \begin{array}{cc}
             0 & B \\
             0 & 0 \\
           \end{array}
         \right)
 \mid B^t=B \in M_2(F) \right\}.
$$
Then $\dim L^{-}_1=1, \dim L^{+}_1=3$.

As a vector space $L$ is embedded into $M_4(F)$. Using ordinary associative 
matrix multiplication we can define super-Lie product on $L$ as
$$
\{x,y\}=xy-(-1)^{|x||y|}yx
$$
for any homogeneous $x,y\in L_0\cup L_1$ where $|x|=0$ if $x\in L_0$ while $|x|=1$, 
if $x\in L_1$.

Note that $L_0$ is a Lie algebra isomorphic to $sl_2(F)$. We will identify $L_0$ 
with $sl_2(F)$ and use the standard basis of $sl_2(F)$
$$
e= \left(
        \begin{array}{cc}
           0 & 1 \\
           0 & 0 \\
    \end{array}
        \right), 
\quad 
f= \left(
                    \begin{array}{cc}
                      0 & 0 \\
                       1 & 0 \\
                      \end{array}
                       \right),
\quad 
h= \left(
            \begin{array}{cc}
              1 & 0 \\
              0 & -1 \\
            \end{array}
          \right).
$$

Furthermore we will not use associative multiplication. This will allow us to omit 
super-Lie brackets, i.e. to write $ab$ instead of $\{a,b\}$. We will also use the 
notation $ab\cdots c$ for the left-normed product $\{\{\ldots\{a,b\},\ldots\},c\}$. 

We will also use the following agreement for denoting alternating sets of variables.
If $f=f(x_1,\ldots,x_k,y_1,\ldots,y_n))$ is a mutilinear polynomial and we apply to $f$ the
operator of alternation on variables $x_1,\ldots, x_k$, then we will write the same symbol
(bar, double bar, tilde, double tilde, etc.) over the variables $x_1,\ldots, x_k$, that is
$$
f(\bar x_1,\ldots,\bar x_k,y_1,\ldots y_n)=\sum_{\sigma\in S_k} ({\rm sgn~\sigma})
f(x_{\sigma(1)},\ldots, x_{\sigma(k)}, y_1,\ldots, y_n).
$$
For example $\bar xa\bar y=xay-yax$, or
$$
{\bar{\bar x}}_1a {\bar{\bar x}}_2 b{\bar{\bar x}}_3=
\sum_{\sigma\in S_3} ({\rm sgn~\sigma})
x_{\sigma(1)} a x_{\sigma(2)} b x_{\sigma(3)}.
$$
We will also use this notation for non-multilinear polynomials with repeating variables
as follows
$$
\bar x_1\bar x_2 a {\bar{\bar x}}_1{\bar{\bar x}}_2= x_1x_2 a {\bar{\bar x}}_1{\bar{\bar x}}_2
-x_2x_1 a _1{\bar{\bar x}}_1{\bar{\bar x}}_2=
$$
$$
x_1x_2ax_1x_2 - x_1x_2ax_2x_1 - x_2x_1ax_1x_2 + x_2x_1ax_2x_1.
$$

Following this agreement we consider an alternating expression 
$h\bar e\bar f\,\bar h$ in the Lie
algebra $L_0=sl_2(F)$. Since $he=2e, hf=-2f, ef=h$, it easily follows that 
$$
h\bar e\bar f\bar h=8h
$$
and
\begin{equation}\label{eq1}
h\underbrace{\bar e\bar f\,\bar h
{\bar{\bar e}}{\bar{\bar f}}\,{\bar{\bar h}}\cdots\tilde e\tilde f
\tilde h }_{t~{\rm alternating~triples}}=8^t h.
\end{equation}

Consider a multilinear polynomial
$$
g=g(x_0,x_1^1,x_2^1,x_3^1,\ldots,x_1^t,x_2^t,x_3^t)=
Alt_1 Alt_2\ldots Alt_t(x_0x_1^1\cdots x_3^1)
$$
where $Alt_j$ is the operator of alternation on $x_1^j,x_2^j,x_3^j$, $1\le j\le t$, that is
$$
g=x_0\bar x_1^1\bar x_2^1\bar x_3^1\cdots \tilde x_1^t\tilde x_2^t\tilde x_3^t.
$$
The evaluation $\varphi:x_1^1,\ldots, x_1^t\to e, x_2^1,\ldots,x_2^t\to f,
x_0,x_3^1,\ldots,x_3^t\to h$ gives us
\begin{equation}\label{eq2}
\varphi(g)=8^th
\end{equation}
in $L_0$ by (\ref{eq1}). Moreover, if we denote  the symmetrization on variables 
$x^1_i,\ldots, x^t_i$, $i=1,2,3$ by $Sym_i$, then it follows  from (\ref{eq2}) and 
the definition of $\varphi$ that
\begin{equation}\label{eq3a}
\varphi(Sym_1Sym_2Sym_3(g))=(t!)^3 8^t h
\end{equation}
in $L_0$ by virtue of (\ref{eq1}). An element
\begin{equation}\label{eq3}
g'=Sym_1Sym_2Sym_3(g)
\end{equation}
with the fixed $x_0$ generates in $P_{3t+1,0}$ an irreducible $S_{3t}$-submodule with the
character $\chi_\lambda,\lambda=(t,t,t)$ where the permutation group $S_{3t}$ acts on
$x_1^1,x_2^1,x_3^1,\ldots,x_1^t,x_2^t,x_3^t$.

Given a partition $\mu=(\mu_1,\ldots, \mu_d) \vdash n$, we define the function
$$
\Phi(\mu)=\frac{1}{z_1^{z_1}\cdots z_d^{z_d}},
$$
where
$$
z_1=\frac{\mu_1}{n}, \ldots, z_d=\frac{\mu_d}{n}.
$$

The value of $\Phi(\mu)$ is closely connected with $\deg\chi(\mu)$.

\begin{lemma}\label{l1}\cite[Lemma 1]{GZ-LMS}
Let $n\ge 100$. Then
$$
\frac{\Phi(\mu)^n}{n^{d^2+d}} \le d_\mu \le n \Phi(\mu)^n.
$$
\end{lemma}
\hfill $\Box$

In particular, if $m=3t$ and $\mu=(t,t,t)$ then
\begin{equation}\label{eq4}
\deg\chi_\mu \ge m^{-12} 3^m.
\end{equation}

In the next step we will construct an irreducible $S_{6k}$-submodule in
$P_{1,6k+1}\not\subset Id^{gr}(L)$ where $S_{6k}$ acts on some
$6k$ odd variables. Denote
$$
d   =  \left(
           \begin{array}{cccc}
             0 & 0 & 0 & 0\\
             0 & 0 & 0 & 0\\
             0 & 1 & 0 & 0\\
             -1 & 0 & 0 & 0\\
           \end{array}
         \right)
\in L_1^-.
$$
Then $L_1^-=Span<d>$. It is not difficult to check that $L_1^+L_1^-=L_0$.
Hence there exist $a,b,c\in L_1^+$ such that
$$
ad=e,\quad bd=f, \quad cd=h
$$
where $\{e,f,h\}$ is the fixed basis of $L_0$. It follows that
$$
h(\bar ad)(\bar bd)(\bar cd)=8h
$$
in $L$ and
$$
h(\bar ad)(\bar bd)(\bar cd) (\bar{\bar a}\bar d)
(\bar{\bar b}d)(\bar{\bar c}d)=
h(\bar ad)(\bar bd)(\bar cd) (\bar{\bar a} d)
(\bar{\bar b}d)(\bar{\bar c}d)=64 h.
$$
Repeating this procedure and using (\ref{eq2}) we obtain a multialternating expression
\begin{equation}\label{eq6}
H=h(\bar ad)(\bar bd)(\bar cd) (\bar{\bar a}\bar d)(\bar{\bar b}d)(\bar{\bar c}d)
(\tilde a{\bar {\bar d}})(\tilde  bd)(\tilde  cd) \cdots
({\tilde{\tilde a}}\hat d)({\tilde{\tilde b}}d)({\tilde{\tilde c}}d)
(a{\tilde{\tilde d}})(bd)(cd)=8^{k+1}h
\end{equation}
depending on one $h$, $k+1$ elements $a,b,c$ and $3(k+1)$ elements $d$.
The element $H$ on the left hand side of (\ref{eq6}) contains  $k$ alternating 
sets $\{a,b,c,d\}$. The first set consists of 1st $a$, 1st $b$, 1st $c$ and
4th $d$. The second set consists of 2nd $a$, 2nd $b$, 2nd $c$ and 7th $d$, and
so on. The element $H$ also contains  $2(k+1)+1$ non-alternating entries $d$
and four extra factors $h,a,b,c$ out of alternating sets. This $H$ is a value
of the following multilinear polynomial: denote by
$$
w=w(x_0,y_1^1,y_2^1,y_3^1,z_1^1,z_2^1,z_3^1,\ldots, y_1^{k+1},y_2^{k+1},y_3^{k+1},
z_1^{k+1},z_2^{k+1},z_3^{k+1})=
$$
$$
Alt_1\ldots Alt_k\left(x_0(y_1^1z_1^1)(y_2^1z_2^1)(y_3^1z_3^1)
(y_1^2z_1^2)(y_2^2z_2^2)(y_3^2z_3^2)\cdots\right.
$$
$$
\left.(y_1^{k+1}z_1^{k+1})(y_2^{k+1}z_2^{k+1})(y_3^{k+1}z_3^{k+1})
\right)
$$
where $x_0$ is an even variable, all $y^\alpha_\beta, z^\alpha_\beta$ are odd
and $Alt_j$ is the operator of alternation on $y_1^j, y_2^j, y_3^j,
z_1^{j+1},j=1,\ldots,k$. Then $\varphi(w)=H$ where $\varphi$ is an evaluation of the form
$$
\varphi(x_0)=h,\,  \varphi(y_1^j)=a,\,  \varphi(y_2^j)=b,\, \varphi(y_3^j)=c,\,
\varphi(z_i^j)=d,\, j=1,\ldots, k+1, i=1,2,3.
$$

Also denote  by $Sym_1,Sym_2,Sym_3$ the symmetrization on the sets
$\{y_1^1,\ldots, y_1^k\}$, $\{y_2^1,\ldots, y_2^k\}$, $\{y_3^1,\ldots, y_3^k\}$,
respectively, and by $Sym_4$ the symmetrization on \\ 
$\{z_1^2,z_2^2,z_3^2\ldots, z_1^{k+1}, z_2^{k+1},z_3^{k+1}\}$.
If
$$
w'=Sym_1 Sym_2 Sym_3Sym_4(w)
$$
then
$$
\varphi(w')=(3k)!(k!)^3\varphi(w)=(3k)!(k!)^3 8^{k+1}h
$$
in $L$. In particular, $w'$ is not an identity of $L$.

Now let the permutation group $S_{6k}$ act on the set
$\{y_i^j, z_i^{j+1}|1\le j \le k, i=1,2,3\}$ and let 
$x_0, z_1^1,z_2^1,z_3^1, y_1^{k+1}, y_2^{k+1}, y_3^{k+1}$ be fixed. Then
$w'$ generates an irreducible $S_{6k}$-submodule in $P_{1, 6(k+1)}$ corresponding
to the partition $\mu=(3k,k,k,k)\vdash r=6k$ with
$$
\Phi(\mu)=2\sqrt 3.
$$
Hence
\begin{equation}\label{eq7}
\deg \chi_\mu \ge r^{-20}(2\sqrt 3)^r
\end{equation}
by Lemma \ref{l1}.

Now let
$$
u=w'(g', y_1^1,y_2^1,y_3^1,\cdots , z_1^{k+1}, z_2^{k+1}, z_3^{k+1})
$$
where $g'$ is taken from (\ref{eq3}). If $m=3t,  r=6k$, group $S_m$ acts on even 
variables $x_i^j$ from $g'$ whereas $S_r$ acts on odd
variables $\{y_i^j,z_i^j\}$ (except $z_1^1,z_2^1,z_3^1, y_1^{k+1}, y_2^{k+1}, y_3^{k+1}$)
then $u$ is not a graded identity of $L$ as follows from (\ref{eq3a}), (\ref{eq4})
and (\ref{eq6}) and it generates an irreducible 
$S_m\times S_r$-submodule $M$ in $P_{m+1,r+6}$ with the character $\chi_{\lambda,\mu}$
where $\lambda=(t,t,t)$, $\mu=(3k,k,k,k)$. Hence by (\ref{eq4}), (\ref{eq7})
$$
\dim M= \deg\chi_\lambda\deg\chi_\mu\ge \frac{1}{(m+r)^{32}}3^m(2\sqrt 3)^r
$$
 and then by (\ref{eq0})
\begin{equation}\label{eq8}
c_{m+1,r+6}(L)=\dim P_{m+1,r+6}(L) \ge \dim M \ge \frac{3^m(2\sqrt 3)^r}{(m+r)^{32}}.
\end{equation}

The inequality (\ref{eq8}) means that we have proved the following lemma.

\begin{lemma}\label{l2}
Let $t,r\ge 1$ be arbitrary integers and $m=3t, r=6k$. Then
$$
c_{m+1,r+6}(L) \ge \frac{3^m(2\sqrt 3)^r}{(m+r)^{32}}.
$$
\end{lemma}
\hfill $\Box$

Now we will find a lower bound for $n$th graded codimension of $L$ for
the special case of $n$.

\begin{lemma}\label{l3}
Let $n-7$  be a multiple of $6$. Then
$$
c_n^{gr}(L) \ge \frac{1}{3^{18}n^{38}}(3+2\sqrt 3)^n.
$$
\end{lemma}
{\em Proof.} Let $q=n-7$. Then applying Lemma \ref{l2} we obtain
$$
c_n^{gr}(L)=\sum_i {n\choose i} c_{i,n-i}(L) \ge\sum_{j=0}^{q/6}
{n\choose 6j+1} c_{1+6j,q+6-6j}(L) \ge 
$$
$$
\frac{1}{q^{32}} \sum_{j=0}^{q/6}
{n\choose 6j+1} 3^{6j}(2\sqrt 3)^{q-6j} \ge \frac{A}{n^{33}}
$$
where
$$
A=\sum_{j=0}^{q/6} {n\choose 6j} 3^{6j}(2\sqrt 3)^{q-6j}
$$
since
$$
{n\choose i} \le n {n\choose i+1}.
$$
Now, since given $0\le j<\frac{q}{6}$, we have
$$
{q\choose 6j+i} 3^{6j+i}(2\sqrt 3)^{q-6j-i} < (3q)^5{q\choose 6j} 3^{6j}(2\sqrt 3)^{q-6j}
$$
for all $1\le i\le 5$. It follows that
$$
A>\frac{6}{(3n)^5}\sum_{i=0}^q {q\choose i} 3^i(2\sqrt 3)^{q-i}=\frac{6}{(3n)^5}
(3+2\sqrt 3)^q >\frac{2}{3^{18}n^5}(3+2\sqrt{3})^n.
$$
Hence
$$
c_n^{gr}(L)>\frac{(3+2\sqrt 3)^n}{3^{18}n^{38}}.
$$
\hfill $\Box$

Now we consider  the case when $n-7$ is not a multiple of $6$, that is 
$n-7\equiv i(mod~6)$ with $1\le i \le 5$.

\begin{lemma}\label{l4}
Let $t,r,i\ge 1$ be arbitrary integers, $m=3t$, $r=6k$ and $i\le 5$. Then
$$
c_{m+1+i,r+6}(L)\ge \frac{3^m(2\sqrt 3)^{r}}{(m+r)^{32}}.
$$
\end{lemma}
{\em Proof.} 
The proof is similar to the proof of Lemma \ref{l2}. We only need to change the
polynomial $u=w(g',y_1^1,\ldots, z_3^{k+1})$ to $u'=ux_1\cdots x_i$ and consider 
an evaluation $\varphi$ with the same values on $x_0,x_j^i, y_j^i, z_j^i$ as in
Lemma \ref{l2} and $\varphi(x_1)=e$, $\varphi(x_2)=\ldots =\varphi(x_i)=h$ if 
$i\ge 2$. Then $\varphi(u')=\pm 2^i\varphi(u)\ne 0$ and we are done.
\hfill $\Box$

Slightly modifying arguments of Lemma \ref{l3} and using Lemma \ref{l4}, 
 we get the following result for arbitrary $1\le i \le 5$.

\begin{lemma}\label{l5}
Let $n-7\equiv i(mod~6), 1\le i \le 5$. Then
$$
c_n^{gr}(L) \ge \frac{1}{3^{18}n^{43}}(3+2\sqrt 3)^n.
$$
\end{lemma}
\hfill $\Box$

Now we are ready to prove the main result of the paper.

\begin{theorem}
Graded PI-exponent of Lie superalgebra $L=b(2)$ exists and is equal to
$$
exp^{gr}(L)=3+2\sqrt 3.
$$
\end{theorem}
{\em Proof.} By Lemmas \ref{l3} and \ref{l5}
$$
\underline{exp}^{gr}(L)=\liminf_{n\to\infty} \sqrt[n]{c_n^{gr}(L)}
\ge 3+2\sqrt 3.
$$
On the other hand, $\overline{exp}^{gr}(b(t))\le t^2-1+t\sqrt{t^2-1}$
for all $t\ge 2$ as proved in \cite{RZ-ART}. Hence the limit
$$
exp^{gr}(b(2))=\lim_{n\to\infty} \sqrt[n]{c_n^{gr}(b(2))}
$$
exists and
$$
exp^{gr}(L)=3+2\sqrt 3.
$$
\hfill $\Box$

\end{document}